\documentclass[12pt]{article}
\begin{document}
\baselineskip=13.5pt
\parskip=8pt
\newcommand{\la}{\langle}
\newcommand{\ra}{\rangle}
\newcommand{\psp}{\vspace{0.4cm}}
\newcommand{\pse}{\vspace{0.2cm}}
\newcommand{\ptl}{\partial}
\newcommand{\dlt}{\delta}
\newcommand{\sgm}{\sigma}
\newcommand{\al}{\alpha}
\newcommand{\be}{\beta}
\newcommand{\G}{\Gamma}
\newcommand{\gm}{\gamma}
\newcommand{\vs}{\varsigma}
\newcommand{\lmd}{\lambda}
\newcommand{\td}{\tilde}
\newcommand{\vf}{\varphi}
\newcommand{\rd}{\mbox{Rad}}
\newcommand{\ad}{\mbox{ad}}
\newcommand{\stl}{\stackrel}
\newcommand{\ol}{\overline}
\newcommand{\ul}{\underline}
\newcommand{\es}{\epsilon}
\newcommand{\dmd}{\diamond}
\newcommand{\clt}{\clubsuit}
\newcommand{\vt}{\vartheta}
\newcommand{\ves}{\varepsilon}
\newcommand{\dg}{\dagger}
\newcommand{\kn}{\mbox{ker}}
\newcommand{\for}{\mbox{for}}
\newcommand{\dvg}{\mbox{div}}
\newcommand{\rar}{\rightarrow}
\newcommand{\NJ}{\mathbb{N}^{\ell_1+\ell_2}}
\newcommand{\ZJ}{\mathbb{Z}^{\ell_1+\ell_2}}
\newcommand{\bs}{\backslash}
\newcommand{\der}{\mbox{Der}\:}
\newcommand{\lra}{\Longleftrightarrow}
\newcommand{\BB}{\mathbb}
\def\a{\alpha}
\def\b{\beta}
\def\N{\hbox{$I\hskip -4pt N$}}
\def\F{\hbox{$I\hskip -4pt F$}}
\def\Z{\hbox{$Z\hskip -5.2pt Z$}}
\def\sZ{\hbox{$\sc Z\hskip -4.2pt Z$}}
\def\Q{\hbox{$Q\hskip -5pt\vrule height 6pt depth 0pt\hskip 6pt$}}
\def\sQ{\hbox{$\sc Q\hskip -4pt\vrule height 6pt depth 0pt\hskip 6pt$}}
\def\R{\hbox{$I\hskip -3pt R$}}
\def\C{\hbox{$C\hskip -5pt \vrule height 6pt depth 0pt \hskip 6pt$}}
\def\J {\rb{5pt}{\mbox{$^{\rar}_{\dis J}$}}}
\def\qed{\hfill \hfill \ifhmode\unskip\nobreak\fi\ifmmode\ifinner
\else\hskip5pt\fi\fi\hbox{\hskip5pt\vrule width4pt height6pt
depth1.5pt\hskip 1 pt}}
\def\d{\delta}
\def\D{\Delta}
\def\g{\gamma}
\def\G{\Gamma}
\def\l{\lambda}
\def\L{\Lambda}
\def\o{\omiga}
\def\p{\psi}
\def\Si{\Sigma}
\def\si{\sigma}
\def\sc{\scriptstyle}
\def\ssc{\scriptscriptstyle}
\def\dis{\displaystyle}
\def\cl{\centerline}
\def\DD{\hbox{$I\hskip -4pt D$}}
\def\sDD{\hbox{$\sc I\hskip -2.5pt D$}}
\def\ll{\leftline}
\def\rl{\rightline}
\def\sF{\hbox{$\sc I\hskip -2.5pt F$}}
\def\nl{\newline}
\def\ol{\overline}
\def\ul{\underline}
\def\wt{\widetilde}
\def\wh{\widehat}
\def\rar{\rightarrow}
\def\Rar{\Rightarrow}
\def\lar{\leftarrow}
\def\Lar{\Leftarrow}
\def\rla{\leftrightarrow}
\def\Rla{\Leftrightarrow}
\def\bs{\backslash}
\def\hs{\hspace*}
\def\vs{\vspace*}
\def\rb{\raisebox}
\def\ra{\rangle}
\def\la{\langle}
\def\SS{\hbox{$S\hskip -6.2pt S$}}
\def\ni{\noindent}
\def\hi{\hangindent}
\def\ha{\hangafter}
\def\Box{\mbox{$|\!\!\!\!|\!\!\!\!|$}}
\def\mathbb#1{{\mbox{c$\!\!\!\!$}\cal #1}}
\def\Bbb#1{{\mbox{b$\!\!\!\!$}\cal #1}}
\def\AA{{\cal A}}
\def\BB{{\cal B}}
\def\ii{{\bf i}}
\def\jj{{\bf j}}
\def\kk{{\bf k}}
\def\mm{{\bf m}}
\def\HOM{\mbox{Hom$^*_{\sZ}(\BB,\F)$}}
\def\sone{{1\hskip -5.5pt 1}}
\def\one{{1\hskip -6.5pt 1}}
\begin{center}{{\large\bf Derivations and structure of the Lie algebras of
Xu type}\footnote{AMS Subject Classification-Primary: 17B40,
17B65, 17B70\nl\hs{4ex}Supported by a grant from Educational
Department of China} }\end{center}
 \vs{-7pt}\par{\small
\cl{\small(appeared in {\it Manuscripta Math.}, {\bf105} (2001),
483-500)} \vs{-3pt}\par \cl{Yucai Su\footnote{Email:
ycsu@sjtu.edu.cn}}\vs{-4pt} \cl{\small Department of Mathematics,
Shanghai Jiaotong University, Shanghai 200030, China}
\par
{\bf Abstract.} We determine the derivation algebras and the
isomorphism classes of a family of the simple Lie algebras
introduced recently by Xu [{\it Manuscripta Math}, {\bf 100},
489--518 (1999)]. The structure space of these algebras is given
explicitly. }\par\ \vs{-19pt}\par\ni {\bf 1. \ Introduction}
\par\ni
Xu [X] recently introduced two new families of infinite-dimensional simple
Lie algebras and a new family of infinite-dimensional simple Lie
superalgebras. These algebras are generalizations of the algebras introduced
by Block [B] and Dokovic and Zhao [DZ].
\par
In this paper, we shall determine the derivation algebras and the
isomorphism classes of the second family of the simple Lie algebras
introduced by Xu. The problems for the first family of the algebras are
settled in [SZ].
\par
Below we shall introduce the normalized form of these algebras. Let $\F$ be
an algebraically closed field with characteristic 0. All the vector spaces
are assumed over \F. Denote by $\Z$ the ring of integers and by $\N$ the
additive semigroup of numbers $\{0,1,2,...\}$.
\par
Take an additive subgroup $\G$ of $\F^4$ such that
$$
\si=(1,0,1,0)\in\G,\ \
\G\cap\F1_{[p]}\ne\{0\},\ \ \G\cap(1_{[q]}+\F1_{[q+1]})\ne\emptyset,
\eqno(1.1)$$
for $p=1,2,3,4,$ $q=1,3$, where here and below, we use the notation
$$
a_{[p]}=(0,...,0,\rb{4pt}{\mbox{$^{^{\sc p}}_{\dis a}$}},0...,0)\in\F^4\ \ \
\for\ \ \ a\in\F,\ p=1,2,3,4.
\eqno(1.2)$$
An element $\a\in\F^4$ will be denoted by $\a=(\a_1,\a_2,\a_3,\a_4)$. For
$p=1,2,3,4$, pick $J_p=\{0\}$ or $\N$, and set
$J=J_1\times J_2\times J_3\times J_4.$ Elements in $J$ will be denoted by
$\ii\!=\!(i_1,i_2,i_3,i_4)$. Denote by $\AA_4\!=\!\AA_4(\G,J)$ the semigroup
algebra of $\G\times J$ with a basis
$\{x^{\a,\ii}\,|\,(\a,\ii)\!\in\!\G\!\times\! J\}$ and the algebraic
operation $\cdot$ defined by:
$$
x^{\a,\ii}\cdot x^{\b,{\ssc\,}\jj}=x^{\a+\b,\ii+\jj}\ \ \for\ \
(\a,\ii),(\b,\jj)\in\G\times J.
\eqno(1.3)$$
Then $x^{0,0}$ is an identity element, which will be simply denoted by 1.
For simplicity, we denote $x^\a=x^{\a,0}$. For $p=1,2,3,4$, we define
derivation $\ptl_p$ of $\AA_4$ by
$$
\ptl_p(x^{\a,\ii})=\a_p x^{\a,\ii}+i_p x^{\a,\ii-1_{[p]}}\ \ \for\ \
(\a,\ii)\in\G\times J,
\eqno(1.4)$$
where we adopt the convention that if a notion is not defined but technically
appears in an expression, we always treat it as zero; for instance,
$x^{\a,-1_{[1]}}=0$ for any $\a\in\G$.
\par
We define the following algebraic operation $[\cdot,\cdot]$ on
$\AA_4=\AA_4(\G,J)$:
$$
[u,v]=x^{\si,0}(\ptl_1(u)\ptl_2(v)-\ptl_1(v)\ptl_2(u))
+(u+\ptl_3(u))\ptl_4(v)-(v+\ptl_3(v))\ptl_4(u),
\eqno(1.5)$$
for $u,v\in\AA_4$. \vs{-7pt}Denote
$$
\si_1=-1_{[3]}=(0,0,-1,0),\ \ \si_2=1_{[1]}-2_{[3]}=(1,0,-2,0).
\vs{-2pt}\eqno(1.6)$$
We treat $x^{\si_i}=0$ if $\si_i\notin\G$ for $i=1,2.$ Then $x^{\si_1,0}$ is
in the center of $\AA_4$. Form a quotient Lie algebra
$\BB_4=\BB_4(\G,J)=\AA_4/\F x^{\si_1}$, whose induced Lie bracket is still
denoted by $[\cdot,\cdot]$.
\par
{\bf Theorem 1.1}. {\it The Lie algebra $(\BB_4,[\cdot,\cdot])$ is simple
if $J\ne\{0\}$ or $\si_2\notin\G$. If $J=\{0\}$ and $\si_2\in\G$, then
$\BB_4^{(1)}=[\BB_4,\BB_4]$ is simple and $\BB_4=\BB_4^{(1)}
\oplus(\F x^{\si_1,0}+\F x^{\si_2,0})$.}
\par
The above theorem was due to Xu [X] and the simple Lie algebras
$\BB_4^{(1)}$, which will now be simply denoted by $\BB$ or $\BB(\G,J)$,
are the normalized form of the second family of the algebras constructed
in [X]. We shall refer these algebras to as the {\it Lie algebras of Xu
type}. We will denote the image of the element $x^{\a,\ii}$ still by
$x^{\a,\ii}$. In particular, we have $x^{\si_1,0}=0$. We set $\G^\#=\G$ if
$J\ne\{0\}$ and $\G^\#=\G\bs\{\si_1,\si_2\}$ otherwise, then $\BB$ has a
basis
$$
B=\{x^{\a,\ii}\,|\,(\a,\ii)\in\G^\#\times J,(\a,\ii)\ne(\si_1,0)\}.
\eqno(1.7)$$
\hs{3ex}
In Section 2, we shall determine the derivations of $\BB$, then in Section
3, we shall give the isomorphism classes and structure space of the Lie
algebras of Xu type. Our main results are Theorems 2.1, 3.2, 3.3.
\par\
\vs{-19pt}\par\ni
{\bf 2. \ Derivations of $\BB$}
\par\ni
We shall determine the structure of the derivation algebra of the Lie
algebra $\BB$.
\par
Recall that a derivation $d$ of the Lie algebra $\BB$ is a linear
transformation on $\BB$ such that $d([u,v])=[d(u),v]+[u,d(v)]$ for
$u,v\in\BB.$ Denote by $\der\BB$ the space of the derivations of $\BB$. It
is well known that $\der\BB$ forms a Lie algebra with respect to the
commutator of linear transformations of $\BB$, and $\ad_{\cal B}$ is an
ideal of $\der\BB$. Elements in $\ad_{\cal B}$ are called {\it inner
derivations}, while elements in $\der\BB\bs\ad_{\cal B}$ are called {\it
outer derivations}.
\par
We rewrite (1.5) as \vs{-2pt}follows:
$$
\matrix{
[x^{\a,\ii},x^{\b,{\ssc\,}\jj}]\!\!\!\!\!&=(\a_1\b_2-\b_1\a_2)x^{\a+\b+\si,\ii+\jj}
+(\a_1j_2-\b_1i_2)x^{\a+\b+\si,\ii+\jj-1_{[2]}}
\vs{4pt}\hfill\cr&
+(i_1\b_2-j_1\a_2)x^{\a+\b+\si,\ii+\jj-1_{[1]}}
+(i_1j_2-j_1i_2)x^{\a+\b+\si,\ii+\jj-1_{[1]}-1_{[2]}}
\vs{4pt}\hfill\cr&
+((\a_3\!+\!1)\b_4\!-\!(\b_3\!+\!1)\a_4)x^{\a+\b,\ii+\jj}
\!+\!((\a_3\!+\!1)j_4\!-\!(\b_3\!+\!1)i_4)x^{\a+\b,\ii+\jj-1_{[4]}}
\vs{4pt}\hfill\cr&
+(i_3\b_4-j_3\a_4)x^{\a+\b,\ii+\jj-1_{[3]}}
+(i_3j_4-j_3i_4)x^{\a+\b,\ii+\jj-1_{[3]}-1_{[4]}},
\hfill\cr
}
\vs{-2pt}\eqno(2.1)$$
for $(\a,\ii),(\b,\jj)\in\G^\#\times J$. For convenience, we fix four
elements in $\G$:
$$
\tau_1=-\si,\ \tau_2=(0,a,-1,0),\ \tau_3=0,\ \tau_4=(0,0,-1,b)\in\G
\mbox{ for some }a,b\in\F\bs\{0\}.
\eqno(2.2)$$
Note that by (1.1) there exist $a',b',c',d',e'\in F$ with $a',b',e'\ne0$
such that
$$
(a',0,0,0),\,(0,b',0,0),\,(1,c',0,0),\,(0,0,1,d'),\,(0,0,0,e')\in\G,
$$
then $(0,kb'+c',-1,0)=k(0,b',0,0)+(1,c',0,0)-\si\in\G$ and
$-(0,0,1,d')+k(0,0,0,e')=(0,0,-1,d'+ke')\in\G$ for all $k\in\Z$; in
particular, by choosing $k$ with $a=kb'+c'\ne0$, $b=d'+ke'\ne0$, it shows
that such $a,b$ exist.
\par
Observing from (1.5), we \vs{-2pt}see
$$
\ad_{x^{\tau_p}}=\left\{\matrix{
-\ptl_2\hfill&\!\!\!\mbox{if \ }p=1,\vs{2pt}\hfill\cr
-a_2x^{\tau_2+\si}\ptl_1
 \hfill&\!\!\!\mbox{if \ }p=2,\vs{2pt}\hfill\cr
\ptl_4\hfill&\!\!\!\mbox{if \ }p=3,\vs{2pt}\hfill\cr
-a_4x^{\tau_4}(\ptl_3+1)\hfill&\!\!\!\mbox{if \ }p=4,\hfill\cr
}\right.
\vs{-2pt}\eqno(2.3)$$
where in general, $x^{\a,\ii}\ptl_p$ is the operator on $\BB:
(x^{\a,\ii}\ptl_p)(x^{\b,{\ssc\,}\jj})=x^{\a,\ii}(\ptl_p(x^{\b,{\ssc\,}\jj}))$,
and $x^{\a,\ii}$ acts on $\BB$ via the multiplication of $\AA_4$.
\par
Denote by $\HOM$ the space of additive functions $\mu:\G\rar\F$ satisfying
$\mu(\si)=0$. For each $\mu\in\HOM$, one can define a linear transformation
$d_\mu$ by
$$
d_\mu(x^{\a,\ii})=\mu(\a)x^{\a,\ii}\ \ \ \for\ \ \ (\a,\ii)\in\G^\#\times J.
\eqno(2.4)$$
Since $\mu(\si)=0$ it is straightforward to verify by (2.1) that
$d_\mu\in\der\BB$. We regard $\HOM$ as a subspace of $\der\BB$ by the
embedding: $\mu\mapsto d_\mu$. Define $\ptl_{t_p}$ by
$$
\ptl_{t_p}(x^{\a,\ii})=i_p x^{\a,\ii-1_{[p]}}\ \ \
\for\ \ \ (\a,\ii)\in\G^\#\times J\mbox{ \ and \ }p=1,2,3,4.
\eqno(2.5)$$
Then by (2.1), one can verify that $\ptl_{t_p}$ are derivations of $\BB$
(of course, $\ptl_{t_p}=0$ if $J_p=\{0\}$). Assume that $\si_1\in\G$.
Observe that for $p=1,2,3,4$, if $J_p=\{0\}$, then $x^{\si_1,1_{[p]}}
\notin\BB$ but $[x^{\si_1,1_{[p]}},\BB]\subset\BB$ and thus we can define
the outer derivations \vs{-2pt}(cf.~(2.1))
$$
d_p=\ad_{x^{\si_1,1_{[p]}}}|_\BB:
x^{\b,{\ssc\,}\jj}\mapsto
\left\{\matrix{
\b_2x^{\b+1_{[1]},{\ssc\,}\jj}+j_2x^{\b+1_{[1]},{\ssc\,}\jj-1_{[2]}}\hfill
  &\mbox{ if }\hfill&p=1,\vs{4pt}\hfill\cr
-\b_1x^{\b+1_{[1]},{\ssc\,}\jj}-j_1x^{\b+1_{[1]},{\ssc\,}\jj-1_{[1]}}\hfill
  &\mbox{ if }\hfill&p=2,\vs{4pt}\hfill\cr
\b_4x^{\b-1_{[3]},{\ssc\,}\jj}+j_4x^{\b-1_{[3]},{\ssc\,}\jj-1_{[4]}}\hfill
  &\mbox{ if }\hfill&p=3,\vs{4pt}\hfill\cr
-(\b_3+1)x^{\b-1_{[3]},{\ssc\,}\jj}-j_3x^{\b-1_{[3]},{\ssc\,}\jj-1_{[3]}}\hfill
 &\mbox{ if }\hfill&p=4,\hfill\cr
}\right.
\vs{-2pt}\eqno(2.6)$$
i.e.,
$d_1=x^{1_{[1]}}\ptl_2,d_2=-x^{1_{[1]}}\ptl_1,d_3=x^{-1_{[3]}}\ptl_4$ and
$d_4=-x^{-1_{[3]}}(\ptl_3+1)$. If $\si_1\notin\G$ or $J_p\ne\{0\}$, we set
$d_p=0$. If $J=\{0\}$ and $\si_2\in\G$, we can define another outer
derivation
$$
\ol d_2=\ad_{x^{\si_2}}|_\BB:
x^{\b}\mapsto \b_2x^{\b+2_{[1]}-1_{[3]}}
-\b_4x^{\b+1_{[1]}-2_{[3]}},
\eqno(2.7)$$
i.e., $\ol d_2=x^{2_{[1]}-1_{[3]}}\ptl_2-x^{1_{[1]}-2_{[3]}}\ptl_4$,
and we set $\ol d_2=0$ if $J\ne\{0\}$ or $\si_2\notin\G$. Then we have
\par\ni
{\bf Theorem 2.1}. {\it $\der\BB$ is spanned by
$\ad_\BB,\HOM,\ptl_{t_p},d_p$ and $\ol d_2$ for $p=1,2,3,4.$}
\par\ni
{\it Proof}.
Let $d\in\der\BB$ and let $D$ be the subspace of $\der\BB$ spanned by the
elements mentioned in the theorem. The proof of the theorem is equivalent
to proving that $d\in D$. We shall prove that after a number steps in each
of which $d$ is replaced by $d-d'$ for some $d'\in D$ the 0 derivation is
obtained and thus proving that $d\in D$. This will be done by a number of
claims.
\par
{\bf Claim 1}. We can suppose $d(x^{-\si})=d(1)=0$.
\par
Assume
$$
d(x^{\tau_p})=\sum_{(\a,\ii)\in M_p}e_{\a,\ii}^{(p)}x^{\a+\tau_p,\ii}
\ \ \for\ \ p=1,2,3,4\ \mbox{ and some }\ e_{\a,\ii}^{(p)}\in\F,
\eqno(2.8)$$
where $M_p=\{(\a,\ii)\in\G\times J\,|\,e^{(p)}_{\a,\ii}\ne0\}$ is a finite
set. We set $e_{\a,\ii}^{(p)}=0$ if $(\a,\ii)\notin M_p$.
\par
First assume that $J_4\ne\{0\}$. By induction on $i_4$, we can write
$x^{\a,\ii}=[1,u]$ for some $u\in\BB$: if $\a_4=0$, then
$x^{\a,\ii}=(i_4+1)^{-1}[1,x^{\a,\ii+1_{[4]}}]$, or otherwise,
$x^{\a,\ii}-\a_4^{-1}[1,x^{\a,\ii}]=\a_4^{-1}i_4x^{\a,\ii-1_{[4]}}$, which
by the inductive assumption can be written as $[1,v]$ for some $v\in\BB$. By
this, we can write $d(1)=[1,u]$ for some $u$, and so replacing $d$ by
$d-\ad_u$ gives $d(1)=0$.
\par
Applying $d$ to $[1,x^{\tau_p}]=0$ for $p=1,2$, we obtain that
$\ptl_4(d(x^{\tau_p}))=0$, i.e.,
$$
\a_4=i_4=0\mbox{ \ \ if \ \ }e_{\a,\ii}^{(p)}\ne0,
\eqno(2.9)$$
for $p=1,2$. Similarly, applying $d$ to $[1,x^{\tau_4}]=bx^{\tau_4}$, we see
that (2.9) also holds for $p=4$. If $J_2\ne\{0\}$, then similarly we can
write $d(x^{-\si})=[x^{-\si},v]$ for some $v$. Since $\ptl_4(d(x^{-\si}))=0$,
we can choose $v$ to satisfy $\ptl_4(v)=0$. Thus if we replace $d$ by
$d-\ad_v$, then we still have $d(1)=0$, and we have $d(x^{-\si})=0$.
\par
So suppose $J_2=\{0\}$. Then each term $x^{\a+\tau_1,\ii}$ with $\a_2\ne0$
appearing on the right-hand side of (2.8) for $p=1$ can be written as
$-\a_2^{-1}[x^{-\si},x^{\a+\tau_1,\ii}]$. Thus again we can replace $d$ by
$d-\ad_v$ for some $v$ which satisfies $\ptl_4(v)=0$, to obtain
$$
\a_2=i_2=0\mbox{ \ \  if \ \ }e_{\a,\ii}^{(p)}\ne0,
\eqno(2.10)$$
holds for $p=1$. Applying $d$ to $[x^{-\si},x^{\tau_2}]=-ax^{\tau_2}$, using
(2.10), we see that $\a_2=i_2=0$ if $e_{\a,\ii}^{(2)}\ne0$, i.e., (2.10)
holds for $p=2$. Applying $d$ to $[x^{-\si},x^{\tau_4}]=0$, we see that
(2.10) holds for $p=4$ also. By calculating the coefficients of
$x^{\a+\tau_2,\ii}$ in $[d(x^{-\si}),x^{\tau_2}]+[x^{-\si},d(x^{\tau_2})]
=-a{\sc\,}d(x^{\tau_2})$, using (2.10), we obtain
$$
(a(\a_1-1)e^{(1)}_{\a,\ii}+a(i_1+1)e^{(1)}_{\a,\ii+1_{[1]}})
-ae^{(2)}_{\a,\ii}=-ae^{(2)}_{\a,\ii},
\eqno(2.11)$$
where, if $J_1=\{0\}$, then the second term vanishes. By induction on $i_1$
ranging from ${\rm max}\{j_1\,|\,(\a,\jj)\in M_1\}$ down to zero, we obtain
$$
e^{(1)}_{\a,\ii}=0\ \mbox{ \ if \ }\ \a_1\ne1\mbox{ \ or \ }i_1\ne0.
\eqno(2.12)$$
By calculating the coefficients of $x^{\a+\tau_4-\si,\ii}$ in
$d([x^{-\si},x^{\tau_4}])=0$, using (2.10), we obtain
$$
b\a_3e^{(1)}_{\a,\ii}+b(i_3+1)e^{(1)}_{\a,\ii+1_{[3]}}=0.
\eqno(2.13)$$
where, the second term is vanishing if $J_3=\{0\}$. By induction on $i_3$
ranging from ${\rm max}\{j_3\,|\,(\a,\jj)\in M_1\}$ down to zero, we obtain
$$
\a_3=i_3=0\ \mbox{ \ if \ }e^{(1)}_{\a,\ii}\ne0.
\eqno(2.14)$$
(2.9), (2.10), (2.12) and (2.14) show that if $e^{(1)}_{\a,\ii}\ne0$, then
$(\a,\ii)=(1_{[1]},0)$, but then $e^{(1)}_{\a,\ii}$ is the coefficient of
the term $x^{1_{[1]}-\si}=x^{\si_1}=0$ (cf (2.8) and the statements before
(1.7)). Thus we have $d(x^{-\si})=0$.
\par
Next assume that $J_4=\{0\}$. If $J_2\ne\{0\}$, by interchanging positions
between $\tau_1$ and $\tau_3$ in the above arguments, we can obtain the
result as above. Thus suppose $J_2=\{0\}$. As the statement before (2.10),
we can assume that $\ptl_2(d(x^{-\si}))=0$. Similarly, we can assume that
$\ptl_4(d(1))=0$. Then as the proof above, we can obtain that (2.9), (2.10)
hold for $p=1,2,3,4$. Then (2.11)-(2.14) also hold and so again we can
assume that $d(x^{-\si})=0$. Similarly, by calculating the coefficients of
$x^{\a+\tau_2,\ii}$ in $d([1,x^{\tau_2}])=0$ and calculating the
coefficients of $x^{\a+\tau_4,\ii}$ in $d([1,x^{\tau_4}])=bx^{\tau_4}$,
we obtain that if $e^{(3)}_{\a,\ii}\ne0$, then $(\a,\ii)=(\si_1,0)$ (which
corresponds to the term $x^{\si_1}=0$), and so we have $d(1)=d(x^{-\si})=0$.
This proves Claim 1.
\par
{\bf Claim 2}. We can assume that $d(x^{\tau_p})=0$ for $p=1,2,3,4.$
\par
By Claim 1, we see that (2.9), (2.10) hold for $p=2,4$. For each term
$x^{\a+\tau_2,\ii}$ appearing in the right-hand side of (2.8), using
$[x^{\tau_2},x^{\a-\si,\ii}]=-a(\a_1-1)x^{\a+\tau_2,\ii}-
ai_1x^{\a+\tau_2,\ii-1_{[1]}}$, we see that one can replace $d$ by
$d-\ad_u$ for some $u$, where $u$ is a linear combination of some
$x^{\a-\si,\ii}$, so that $u$ satisfies $\ptl_2(u)=\ptl_4(u)=0$, to obtain
(here we would like to remark that in case $J=\{0\},\a=2_{[1]}-1_{[3]}$,
then $\a-\si=\si_2\notin\G^\#$, but $\ad_{x^{\a-\si}}=\ol d_2$ by (2.7),
this is where we make use of $\ol d_2$),
$$
\a_1=1,\,i_1=0\mbox{ \ \  if \ \ }e_{\a,\ii}^{(2)}\ne0.
\eqno(2.15)$$
By calculating the coefficients of $x^{\a+\tau_2+\tau_4+\si,\ii}$ in
$d([x^{\tau_2},x^{\tau_4}])=0,$ using (2.15), we obtain
$$
\a_1=i_1=0\mbox{ \ \  if \ \ }e_{\a,\ii}^{(4)}\ne0.
\eqno(2.16)$$
Then, using $[x^{\tau_4},x^{\a,\ii}]=-b(\a_3+1)x^{\a+\tau_4,\ii}
-bi_3x^{\a+\tau_4,\ii-1_{[3]}}$, we see that one can replace $d$ by
$d-\ad_u$ for some $u$ which satisfies $\ptl_1(u)=\ptl_2(u)=\ptl_4(u)=0$,
to obtain
$$
\a_3=-1,\,i_3=0\mbox{ \ \  if \ \ }e_{\a,\ii}^{(4)}\ne0.
\eqno(2.17)$$
Then by $d([x^{\tau_2},x^{\tau_4}])=0$, it gives
$$
\a_3=i_3=0\mbox{ \ \  if \ \ }e_{\a,\ii}^{(2)}\ne0.
\eqno(2.18)$$
Using (2.15)-(2.18), we see that $M_2$ or $M_4$ is either empty, or a
singleton $\{(1_{[1]},0)\}$ or $\{(\si_1,0)\}$, thus by replacing $d$ by
$$
d-a^{-1}e^{(2)}_{1_{[1]},0}{\sc\,}\ad_{x^{\si_1,1_{[1]}}}
-b^{-1}e^{(4)}_{\si_1,0}{\sc\,}\ad_{x^{\si_1,1_{[3]}}},
\eqno(2.19)$$
where, the last two terms are either in $\ad_\BB$ or in $\F d_1+\F d_3$
(cf.~(2.6)), we obtain that $d(x^{\tau_p})=0$ for $p=1,2,3,4$. This proves
Claim 2.
\par
Now for any $\b\in\G\bs\{\si_1,\si_2\}$, assume
$$
d(x^\b)=\sum_{(\a,\ii)\in M_\b}c^{(\b)}_{\a,\ii}x^{\a+\b,\ii}
\ \mbox{ for some }\ c^{(\b)}_{\a,\ii}\in\F,
\eqno(2.20)$$
where $M_\b=\{(\a,\ii)\in\G\times J\,|\,c^{(\b)}_{\a,\ii}\ne0\}$ is a finite
set. Again we set $c_{\a,\ii}^{(\b)}=0$ if $(\a,\ii)\notin M_\b$. By
applying $d$ to $[x^{-\si},x^\b]=-\b_2x^\b,[1,x^\b]=\b_4x^\b$, by Claim 1,
we obtain
$$
\a_2=\a_4=i_2=i_4=0\mbox{ \ \  if \ \ }c_{\a,\ii}^{(\b)}\ne0.
\eqno(2.21)$$
Furthermore, if $\b\in\G$ with $\b_1=0$, then by applying $d$ to
$[x^{\tau_2},x^\b]=0$, we obtain
$$
\a_1=i_1=0\mbox{ \ \  if \ \ }c_{\a,\ii}^{(\b)}\ne0\ \mbox{ and }\
\b\in\G\mbox{ \ with \ }\b_1=0.
\eqno(2.22)$$
By (1.1), we can fix an element $\g=(0,0,c,0)\in\G$ with
$c\in\F\bs\{0,-1,-2\}$.
\par
{\bf Claim 3}. We can assume that $d(x^\g)=0$.
\par
If $c_{0,0}^{(\g)}\ne0$, we can define an additive function $\mu\in\HOM$
satisfying that $\mu(\tau_p)=0$ for $p=1,2,3,4$ (cf.~(2.2)) and that
$\mu(\g)=c_{0,0}^{(\g)}$, as follows:
$$
\mu:\a\mapsto c^{-1}c_{0,0}^{(\g)}(-\a_1+a^{-1}\a_2+\a_3+b^{-1}\a_4),
\eqno(2.23)$$
so if we replace $d$ by $d-d_\mu$, $c_{0,0}^{(\g)}$ is vanishing. Since
$d_\mu(x^{\tau_p})=0$, we still have Claim 2 after this replacement.
Furthermore, if $\si_1\in\G$ and $c_{\si_1,0}^{(\g)}\ne0$, then we can
replace $d$ by $d+c_{\si_1,0}^{(\g)}(c+1)^{-1}\ad_{x^{\si_1,1_{[4]}}}$ so
that $c_{\si_1,0}^{(\g)}$ vanishes also. Noting that $[x^{\si_1,1_{[4]}},1]=
-x^{\si_1}=0$ according to the statements before (1.7) and that
$[x^{\si_1,1_{[4]}},x^{\tau_p}]=0$ for $p=1,2,4$, we still have Claim 2 after
this replacement. Therefore we can suppose
$$
c_{0,0}^{(\g)}=c_{\si_1,0}^{(\g)}=0.
\eqno(2.24)$$
Let $j,k\in\Z$. Applying $d$ to $[x^{\g+j\tau_4},x^{k\tau_4}]=
b((c+1)k-j)x^{\g+(j+k)\tau_4}$, where, the coefficient
$b((c+1)k-j)=(c-j+1)kb-(-k+1)jb$ is obtained from the 5th term of the
right-hand side of (2.1), and by calculating the coefficients of
$x^{\a+\g+(j+k)\tau_4,\ii}$, using (2.21) and canceling the common factor
$b$, we \vs{-2pt}obtain
$$
\matrix{
((c+1)k-j)c_{\a,\ii}^{(\g+(j+k)\tau_4)}=
\!\!\!\!&((\a_3+c+1)k-j)c_{\a,\ii}^{(\g+j\tau_4)}
+(i_3+1)kc_{\a,\ii+1_{[3]}}^{(\g+j\tau_4)}
\vs{4pt}\hfill\cr&
+((c+1)k-j(\a_3+1))c_{\a,\ii}^{(k\tau_4)}
-j(i_3+1)c_{\a,\ii+1_{[3]}}^{(k\tau_4)},
\hfill\cr}
\vs{-2pt}\eqno(2.25)$$
where, if $J_3=\{0\}$, the second and the last terms of the right-hand side
of (2.25) are zero. Assume that $c_{\a,\ii}^{(\g)}\ne0$ for some $(\a,\ii)$.
Since we shall be only involved in the choices of $j,k$ with
$-3\le j,k,j+k\le3$, we can set
$$
\mm={\rm max}\,\{\ii\ |\ c_{\a,\ii}^{(\g+j\tau_4)}\ne0
\mbox{ \ or \ }c_{\a,\ii}^{(j\tau_4)}\ne0
\mbox{ \ for some \ }j=-3,...,3\}.
\eqno(2.26)$$
Taking $j=0,\ii=\mm$ in (2.25), and canceling the common factor $k$, we have
$$
(c+1)c_{\a,\mm}^{(\g+k\tau_4)}=
(\a_3+c+1)c_{\a,\mm}^{(\g)}+(c+1)c_{\a,\mm}^{(k\tau_4)}\mbox{ \ if \ }k\ne0.
\eqno(2.27)$$
Multiplying (2.25) by $c+1$ and using (2.27) to substitute the left-hand
side and the first term of the right-hand side, we \vs{-2pt}obtain
$$
\matrix{
k\a_3(\a_3+c+1)c_{\a,\mm}^{(\g)}
=\!\!\!\!&(c+1)(((c+1)k-j)c_{\a,\mm}^{((j+k)\tau_4)}
\vs{4pt}\hfill\cr&
-((\a_3+c+1)k-j)c_{\a,\mm}^{(j\tau_4)}
-((c+1)k-j(\a_3+1))c_{\a,\mm}^{(k\tau_4)}),
\hfill\cr}
\vs{-2pt}\eqno(2.28)$$
if $j,k+j\ne0$. Interchanging $j$ and $k$, and adding the result to (2.28),
canceling the common factor $j+k$, we obtain
$$
\a_3(\a_3+c+1)c_{\a,\mm}^{(\g)}=c(c+1)(c_{\a,\mm}^{((j+k)\tau_4)}
-c_{\a,\mm}^{(j\tau_4)}-c_{\a,\mm}^{(k\tau_4)})\mbox{ \ if \ }j,k,j+k\ne0.
\eqno(2.29)$$
Denote by $c^*$ the left-hand side factored by $c(c+1)$, which does not
depend on $j,k$. Noting that by Claim 2, $c_{\a,\mm}^{(\tau_4)}=0$, using
this in (2.29), we can solve
$$
\a_3(\a_3+c+1)c_{\a,\mm}^{(\g)}=c(c+1)c^*,\ \
c_{\a,\mm}^{(k\tau_4)}=(k-1)c^*\ \ \for\ \ k\ne0.
\eqno(2.30)$$
Using (2.30) in (2.28), by letting $j\ne k$, we obtain
$$
(\a_3+1)c^*=0.
\eqno(2.31)$$
\hs{3ex}
Assume that $\a_3(\a_3+1)(\a_3+c+1)\ne0$. (2.30), (2.31) show that
$c_{\a,\mm}^{(\g)}=c_{\a,\mm}^{(k\tau_4)}=0$,
and that $c_{\a,\mm}^{(\g+k\tau_4)}=0$ by (2.27), this is contrary to (2.26).
Thus this case dos not occur.
\par
Next assume that $\a_3=0$. Then (2.21), (2.22) show that $\a=0$. Then (2.30)
gives
$$
c_{\a,\mm}^{(k\tau_4)}=0.
\eqno(2.32)$$
If $\mm=0$, then the assumption that $c_{\a,\ii}^{(\g)}\ne0$ contradicts
(2.24). So $\mm\ne0$. Then (2.21), (2.22) mean that $m_3\ne0$ and so
$J_3\ne\{0\}$. By (2.32), setting $\ii=\mm$ in (2.25), we obtain that
$((c+1)k-j)c_{\a,\mm}^{(\g+(j+k)\tau_4)}=
((c+1)k-j)c_{\a,\mm}^{(\g+j\tau_4)}.$ Setting $j=0$, this gives
$$
c_{\a,\mm}^{(\g+k\tau_4)}=c_{\a,\mm}^{(\g)}.
\eqno(2.33)$$
Thus the left-hand side does not depend on $k$. Taking $\ii=\mm-1_{[3]}$ in
(2.25), using (2.32), (2.33), we obtain
$$
((c+1)k-j)(c_{\a,\mm-1_{[3]}}^{(\g+(j+k)\tau_4)}
-c_{\a,\mm-1_{[3]}}^{(\g+j\tau_4)}-c_{\a,\mm-1_{[3]}}^{(k\tau_4)})
=m_3kc_{\a,\mm}^{(\g)}.
\eqno(2.34)$$
Setting $k=1$, noting that by Claim 2, $c_{\a,\mm-1_{[3]}}^{(\tau_4)}=0$, we
\vs{-2pt}obtain
$$
\matrix{
(c+1-j)(c_{\a,\mm-1_{[3]}}^{(\g+(j+1)\tau_4)}
-c_{\a,\mm-1_{[3]}}^{(\g+j\tau_4)})=m_3c_{\a,\mm}^{(\g)},
\vs{4pt}\hfill\cr
(c-j)(c_{\a,\mm-1_{[3]}}^{(\g+(j+2)\tau_4)}
-c_{\a,\mm-1_{[3]}}^{(\g+(j+1)\tau_4)})=m_3c_{\a,\mm}^{(\g)},
\hfill\cr}
\eqno\matrix{(2.35)\vs{4pt}\cr(2.36)\cr}\vs{-2pt}$$
where (2.36) is obtained from (2.35) by replacing $j$ by $j+1$. Multiplying
(2.36) by $c+1-j$, multiplying (2.35) by $c-j$, adding the results to kill
the term $c_{\a,\mm-1_{[3]}}^{(\g+(j+1)\tau_4)}$, we obtain
$$
(c-j)(c+1-j)(c_{\a,\mm-1_{[3]}}^{(\g+(j+2)\tau_4)}
-c_{\a,\mm-1_{[3]}}^{(\g+j\tau_4)})=m_3(2c+1-2j) c_{\a,\mm}^{(\g)}.
\eqno(2.37)$$
Setting $k=2$ in (2.34), multiplying it by $(c-j)(c+1-j)$, multiplying
(2.37) by $2(c+1)-j$, subtracting the obtained results from each other, we
obtain
$$
-(c-j)(c+1-j)(2(c+1)-j)c_{\a,\mm-1_{[3]}}^{(2\tau_4)}
=m_3((2c+3)j-2(c+1)^2)c_{\a,\mm}^{(\g)}.
\eqno(2.38)$$
The left-hand side is a polynomial on $j$ of degree 3, while the right-hand
side is of degree 1. This forces $c_{\a,\mm}^{(\g)}=0$. This, together with
(2.32), (2.33), contradicts (2.26). Thus this case dos not occur either.
\par
Next assume that $\a_3=-1$. Then (2.21), (2.22) show that $\a=\si_1$. If
$\mm=0$, then (2.24) contradicts the assumption that $c_{\a,\ii}^{(\g)}\ne0$.
So $\mm\ne0$ and then $J_3\ne\{0\}$. Using (2.30) in (2.27), we obtain
$$
c_{\a,\mm}^{(\g+k\tau_4)}=(k-c-1)c^*.
\eqno(2.39)$$
Taking $\ii=\mm-1_{[3]}$ in (2.25), using (2.30), (2.39), we obtain
$$
((c\!+\!1)k\!-\!j)c_{\a,\mm-1_{[3]}}^{(\g+(j+k)\tau_4)}
\!-\!(ck\!-\!j)c_{\a,\mm-1_{[3]}}^{(\g+j\tau_4)}
\!-\!(c\!+\!1)kc_{\a,\mm-1_{[3]}}^{(k\tau_4)}
\!=\!m_3(\!j\!-\!k(c\!+\!1))c^*,
\eqno(2.40)$$
for $k\ne0$. Setting $k=1$, we \vs{-2pt}obtain
$$
\matrix{
(c+1-j)c_{\a,\mm-1_{[3]}}^{(\g+(j+1)\tau_4)}
-(c-j)c_{\a,\mm-1_{[3]}}^{(\g+j\tau_4)}
=m_3(j-(c+1))c^*,
\vs{4pt}\hfill\cr
(c-j)c_{\a,\mm-1_{[3]}}^{(\g+(j+2)\tau_4)}
-(c-j-1)c_{\a,\mm-1_{[3]}}^{(\g+(j+1)\tau_4)}
=m_3(j-c)c^*,
\hfill\cr
}
\eqno\matrix{(2.41)\vs{4pt}\cr(2.42)\cr}\vs{-2pt}$$
where (2.42) is obtained from (2.41) by replacing $j$ by $j+1$. Multiplying
(2.41) by $c-j-1$ and (2.42) by $c+1-j$ and adding the results together, we
obtain
$$
(c+1-j)(c-j)c_{\a,\mm-1_{[3]}}^{(\g+(j+2)\tau_4)}
-(c-j-1)(c-j)c_{\a,\mm-1_{[3]}}^{(\g+j\tau_4)}=m_3(-2(j-c)^2+j-c+1)c^*.
\eqno(2.43)$$
Setting $k=2$ in (2.40), we have
$$
(2(c+1)-j)c_{\a,\mm-1_{[3]}}^{(\g+(j+2)\tau_4)}
-(2c-j)c_{\a,\mm-1_{[3]}}^{(\g+j\tau_4)}
-2(c+1)c_{\a,\mm-1_{[3]}}^{(2\tau_4)}
=m_3(j-2(c+1))c^*.
\eqno(2.44)$$
Multiplying (2.43) by $-(2(c+1)-j)$ and (2.44) by $(c+1-j)(c-j)$ and adding
the results together, we \vs{-2pt}obtain
$$
\matrix{
(c-j)(c_{\a,\mm-1_{[3]}}^{(\g+j\tau_4)}
+(c+1-j)c_{\a,\mm-1_{[3]}}^{(2\tau_4)})
\vs{4pt}\hfill\cr
=-m_3(j-2(c+1))(j-c+1)(j-c-1)(2(c+1))^{-1}c^*.
\hfill\cr
}
\vs{-2pt}\eqno(2.45)$$
Multiplying (2.41) by $c-j-1$ and using (2.45) to substitute both terms of
the left-hand side of the resulting expression, noting that the coefficient
of $c_{\a,\mm-1_{[3]}}^{(2\tau_4)}$ becomes vanishing, we \vs{-2pt}obtain
$$
\matrix{
m_3(j\!-\!c\!-\!1)((j\!-\!2c\!-\!1))(j\!-\!c\!+\!2)(j\!-\!c)
\vs{4pt}\hfill\cr
\ \ \ \ \ \ \ \ \ \ \
+(c\!-\!j\!-\!1)(j\!-\!2c\!-\!2)(j\!-\!c\!+\!1))(2(c\!+\!1))^{-1}c^*
=m_3(j\!-\!c\!-\!1)(c\!-\!j\!-\!1)c^*.
\hfill\cr}
\vs{-2pt}\eqno(2.46)$$
This shows that $c^*=0$, and so we again get a contradiction to (2.26).
\par
Finally assume that $\a_3=-c-1$. Then the first equation of (2.30) shows
that $c^*=0$, the second shows that $c_{\a,\mm}^{(k\tau_4)}=0$, and (2.27)
shows that $c_{\a,\mm}^{(\g+k\tau_4)}=0$ for $k\ne0$, and then taking
$k=-j\ne0$ in (2.25) gives that $c_{\a,\mm}^{(\g)}=0$, this is again
contrary to (2.26). This shows that it is wrong to make the assumption that
$c_{\a,\ii}^{(\g)}\ne0$ for some $(\a,\ii)$. Thus $d(x^\g)=0$. This proves
Claim 3.
\vs{-1pt}\par
By the statement after (2.2), we can take $e=ka'-1\ne0,\pm1$ for some
$k\in\Z$ such that $\eta=(e,0,-1,0)=k(a',0,0,0)-\si\in\G$. Fix such an
$\eta$.
\par
Applying $d$ to $[x^{\tau_i},x^\eta]=[x^\g,x^\eta]=0,i=1,3,4$, we obtain
$(\a_2,\a_3,\a_4)=(i_2,i_3,i_4)=0$ if $c_{\a,\ii}^{(\eta)}\ne0$. If
$\si_1\in\G$ and $c_{\si+\si_1,0}^{(\eta)}\ne0$, then we can replace $d$ by
$d+c_{\si+\si_1,0}^{(\eta)}e^{-1}\ad_{x^{\si_1,1_{[2]}}}$ (this does not
affect Claims 1,2,3) so that $c_{\si+\si_1,0}^{(\eta)}$ becomes zero. Now
using similar arguments to those given after (2.24) in the proof of Claim 3,
we can prove
$$
\a=\ii=0\mbox{ \ if \ }c_{\a,\ii}^{(\eta)}\ne0.
\eqno(2.47)$$
In fact, if we use the new symbol $y^{\a,\ii}$ to denote $x^{-\a-\si,\ii}$,
then in case $\ii=\jj=0$, (2.1) becomes
$[y^{\a},y^{\b}]=((\a_1+1)\b_2-(\b_1+1)\a_2)y^{\a+\b}+
(\a_3\b_4-\b_3\a_4)y^{\a+\b+\si}$. This shows that there is a symmetry
between $\eta$ and $\g$.
\par
{\bf Claim 4}. We can suppose $d(x^\b)=0$ for all
$\b\in\G\bs\{\si_1,\si_2\}$.
\par
Suppose $\b_4\ne0$. Applying $d$ to
$[x^{\tau_4},x^\b]=-b(\b_3+1)x^{\b+\tau_4},
[x^{\g},x^\b]=(c+1)\b_4x^{\b+\g},$ using Claim 2, Claim 3, (2.21), canceling
the common factors $-b$ and $(c+1)\b_4$, we obtain \vs{-2pt}respectively
$$
\matrix{
\hfill(\b_3+\a_3+1)c_{\a,\ii}^{(\b)}+(i_3+1)c_{\a,\ii+1_{[3]}}^{(\b)}
\!\!\!&=&\!\!\!(\b_3+1)c_{\a,\ii}^{(\b+\tau_4)},
\vs{4pt}\hfill\cr
\hfill c_{\a,\ii}^{(\b)}\!\!\!&=&\!\!\!c_{\a,\ii}^{(\b+\g)}.
\hfill\cr
}
\eqno\matrix{(2.48)\vs{4pt}\cr(2.49)\cr}\vs{-2pt}$$
Replacing $\b$ by $\b+\g$ in (2.48), using (2.49), we have
$$
(c+\b_3+\a_3+1)c_{\a,\ii}^{(\b)}
+(i_3+1)c_{\a,\ii+1_{[3]}}^{(\b)}
=(c+\b_3+1)c_{\a,\ii}^{(\b+\tau_4)}.
\eqno(2.50)$$
This together with (2.48) gives $c_{\a,\ii}^{(\b+\tau_4)}=c_{\a,\ii}^{(\b)}$,
and then (2.48) becomes
$$
\a_3c_{\a,\ii}^{(\b)}+(i_3+1)c_{\a,\ii+1_{[3]}}^{(\b)}=0.
\eqno(2.51)$$
As in (2.13), this shows
$$
\a_3=i_3=0\mbox{ \ if \ }c_{\a,\ii}^{(\b)}\ne0.
\eqno(2.52)$$
If $\b_4=0$, replacing $\b$ by $\b-\tau_4$ in the arguments above, we still
have (2.52). Similarly, using $\tau_2,\eta$ in placing of $\tau_4,\g$ in the
above arguments, we can prove that $\a_1=i_1=0$ if $c_{\a,\ii}^{(\b)}\ne0$.
This and (2.52), (2.21) show that $M_\b$ is either empty or a singleton
$\{(0,0)\}$. Therefore, we can rewrite (2.20) as $d(x^\b)=\mu(\b)x^\b$ by
letting $\mu(\b)=c_{0,0}^{(\b)}$. Using this, by setting $\ii=\jj=0$ in (2.1)
and applying $d$ to it, we \vs{-2pt}obtain
$$
\matrix{
\hfill(\a_1\b_2-\b_1\a_2)(\mu(\a)+\mu(\b))
\!\!\!&=&\!\!\!(\a_1\b_2-\b_1\a_2)\mu(\a+\b+\si),
\vs{4pt}\hfill\cr
\hfill((\a_3+1)\b_4-(\b_3+1)\a_4)(\mu(\a)+\mu(\b))\!\!\!&=&\!\!\!
((\a_3+1)\b_4-(\b_3+1)\a_4)\mu(\a+\b).
\hfill\cr}
\vs{-2pt}\eqno(2.53)$$
From this, one can easily obtain that $\mu$ is an additive function such
that $\mu(\si)=0$. Thus $\mu\in\HOM$, and by replacing $d$ by $d-d_\mu$, we
have Claim 4. This also proves the theorem if $J=\{0\}$.
\par
{\bf Claim 5}. If $J\ne\{0\}$, then $d\in\sum_{p=1}^4\F\ptl_{t_p}$.
\par
Say, $J_1\ne\{0\}$. Assume
$$
d(x^{\b,1_{[1]}})
=\sum_{(\a,\ii)\in M_{\b,1}}c^{(\b,1)}_{\a,\ii}x^{\a+\b,\ii}
\ \ \for\ \ \b\in\G\ \mbox{ and some }\ c^{(\b,1)}_{\a,\ii}\in\F,
\eqno(2.54)$$
where $M_{\b,1}=\{(\a,\ii)\in\G\times J\,|\,c^{(\b,1)}_{\a,\ii}\ne0\}$
is a finite set. Again we set $c^{(\b,1)}_{\a,\ii}=0$ if $(\a,\ii)\notin
M_{\b,1}$. Denote $\BB^{(0)}=
{\rm span}\{x^\a\,|\,\a\in\G\bs\{\si_1,\si_2\}\}.$ Then by Claim 4,
$d(\BB^{(0)})=0$. Applying $d$ to
$[x^{-\si},x^{\tau_2,1_{[1]}}]=-ax^{\tau_2,1_{[1]}}$ and
$[x^{\tau_p},x^{\tau_2,1_{[1]}}]\in\BB^{(0)}$ for $p=2,3,4$, we obtain as
before (for instance, (2.13)) that $M_{\tau_2,1}=\emptyset$ or a singleton
$\{(0,0)\}$. Replacing $d$ by $d-c_{0,0}^{(\tau_2,1)}\ptl_{t_1}$, we can
suppose $d(x^{\tau_2,1_{[1]}})=0$. Similarly, if $J_p\ne\{0\}$ for
$p=2,3,4$, by replacing $d$ by $d-f\ptl_{t_p}$ for some $f\in\F$, we can
suppose $d(x^{\eta,1_{[2]}})=0$ if $p=2$, or $d(x^{\tau_4,1_{[3]}})=0$ if
$p=3$, or $d(x^{\g,1_{[4]}})=0$ if $p=4$. Note that $\BB$ is generated by
$\BB^{(0)}\cup(\BB\cap\{x^{\tau_2,1_{[1]}},x^{\eta,1_{[2]}},
x^{\tau_4,1_{[3]}},x^{\g,1_{[4]}}\})$. Since a derivation is determined
by its actions on generators, we obtain $d(\BB)=0$. This proves Claim 5 and
the theorem.
\qed\par
\par\
\vs{-19pt}\par\ni
{\bf 3. \ Structure of $\BB$}
\par\ni
Before stating the main result of this paper, we need one more lemma. First
recall that a linear transformation $T$ on a vector space $V$ is {\it
locally finite} if for any $u\in V$, the subspace
$$
V_u={\rm span}\{T^n(u)\,|\,n\in\N\},
\eqno(3.1)$$
is finite dimensional; {\it locally nilpotent} if for any $u\in V$ there
exists $n\in\N$ such that $T^n(u)=0$; {\it semi-simple} if for any
$u\in V$, $V_u$ has a basis consisting of eigenvectors of $T$.
\par\ni
{\bf Lemma 3.1}. {\it The sets of locally nilpotent, semi-simple and locally
finite elements of $\der\BB$ are respectively
$$
(\der\BB)_{\it n}=\sum_{p=1}^4\F\ptl_{t_p},\
(\der\BB)_{\it s}=\HOM,\
(\der\BB)_{\it f}=(\der\BB)_{\it n}+(\der\BB)_{\it s}.
\eqno(3.2)$$
Note that $\ad_{x^{-\si}},\ad_1\in\HOM+\sum_{p=1}^4\F\ptl_{t_p}$.
}\par\ni
{\it Proof}.
Note that a linear transformation on the vector space $\BB$ of the form
$$
T=\sum c_{\a,\ii,p}x^{\a,\ii}\ptl_p+
\sum c'_{\a,\ii}x^{\a,\ii}+\sum c_p\ptl_{t_p},
\eqno(3.3)$$
where $x^{\a,\ii}$ acts on $\BB$ via the multiplication of $\AA_4$, is
locally finite if and only if
$$
c_{\a,\ii,p}=c'_{\a,\ii}=0\ \ \mbox{for all}\ \ (\a,\ii)\ne(0,0).
\eqno(3.4)$$
See, for example, the proof of [SX, Lemma 4.1]. Observe from (2.1)
\vs{-2pt}that
$$
\matrix{
\ad_{x^{\a,\ii}}=
\!\!\!\!&x^{\a+\si,\ii}(\a_1\ptl_2-\a_2\ptl_1)
+i_1x^{\a+\si,\ii-1_{[1]}}\ptl_2-i_2x^{\a+\si,\ii-1_{[2]}}\ptl_1
\vs{4pt}\hfill\cr&
+x^{\a,\ii}((\a_3+1)\ptl_4-\a_4(\ptl_3+1))+i_3x^{\a,\ii-1_{[3]}}\ptl_4
-i_4x^{\a,\ii-1_{[3]}}(\ptl_3+1).
\hfill\cr}
\vs{-2pt}\eqno(3.5)$$
From this, one immediately sees that if
$d=\sum c_{\a,\ii}{\ssc\,}\ad_{x^{\a,\ii}}+\sum c_p\ptl_{t_p}+d_\mu$, where
$c_{\a,\ii},c_p\in F$ and $d_\mu$ is defined in (2.4) with $\mu\in\HOM$
(recall that $d_1,d_2,d_3,d_4,\ol d_2$ has the form $\ad_{x^{\a,\ii}}$,
$\sc\,$cf.~(2.6),~(2.7)$\sc\,$) is locally finite if and only if
$$
c_{\a,\ii}=0\mbox{ \ for all \ }(\a,\ii)\ne(0,0),(-\si,0).
\eqno(3.6)$$
From this, one can easily obtain the lemma.
\qed\par\ni
{\bf Theorem 3.2}. {\it The Lie algebras $\BB=\BB(\G,J)$ and
$\BB'=\BB(\G',J')$ are isomorphic if and only if there exist
$a_1,a_2,a_3,a_4\in\F$ with $a_2,a_4\ne0$ and $a_1=0$ if $J_1=\{0\}\ne J_2$
and $a_3=0$ if $J_3=\{0\}\ne J_4$, such that either (i) $J=J'$ and the map
$$\phi:
\b\!=\!(\b_1,\b_2,\b_3,\b_4)\mapsto\b'\!=\!
\b\pmatrix{
^{\dis1\ \ 0}_{\dis a_1\ a_2}
\!\!\!\!&0\cr0\!\!\!\!&
^{\dis1\ \ 0}_{\dis a_3\ a_4}
\cr}
\!=\!(\b_1\!+\!a_1\b_2,a_2\b_2,\b_3\!+\!a_3\b_4,a_4\b_4),
\eqno(3.7)$$
is a group isomorphism $\G\cong\G'$; or else (ii)
$(J_1,J_2,J_3,J_4)=(J'_3,J'_4,J'_1,J'_2)$
and the map
$$\phi:
\b\!=\!(\b_1,\b_2,\b_3,\b_4)\mapsto\b'\!=\!
\b\pmatrix{0\!\!\!\!&
^{\dis1\ \ 0}_{\dis a_1\ a_2}
\cr
^{\dis1\ \ 0}_{\dis a_3\ a_4}
\!\!\!\!&0\cr}
\!=\!(\b_3\!+\!a_3\b_4,a_4\b_4,\b_1\!+\!a_1\b_2,a_2\b_2),
\eqno(3.8)$$
is a group isomorphism $\G\cong\G'$.
}
\par\ni
{\it Proof}.
We shall ALWAYS use the symbol $y$ in place of $x$ for the algebra $\BB'$
and use the same notation with a prime to denote any other element
associated with $\BB'$.
\par
``$\Lar$'':
For convenience, we denote $t_p=x^{0,1_{[p]}}$ for $p=1,2,3,4$. Recall
(1.3), we have the commutative associative algebra structure $(\AA_4,\cdot)$.
We define two algebra structures $(\AA_4,\odot_1)$ and $(\AA_4,\odot_2)$
(which are not necessarily associative) by $u\odot_1v=
x^\si\ptl_1(u)\ptl_2(v)$ and $u\odot_2v=(\ptl_3(u)+u)\ptl_4(v)$.
Then we have $[u,v]=u\odot_1v+u\odot_2v-v\odot_1u-v\odot_2u$ and
$$
x^{\a,\ii}\odot_1 x^{\b,{\ssc\,}\jj}=
x^{\a+\b+\si,\ii+\jj-1_{[1]}-1_{[2]}}(\a_1t_1+i_1)(\b_2t_2+j_2)
\mbox{ \ if \ }i_1+j_1,i_2+j_2\ge1.
\eqno(3.9)$$
\par
Assume that we have case (i). Define a character of $\G$, i.e., a
multiplicative function $\chi:\G\rar\F\bs\{0\}$ such that
$\chi(\si)=a_2a_4^{-1}$. We prove that such a character exists: suppose
$\D$ is the maximal subgroup of $\G$ containing $\si$ such that such
character exists for $\D$. If $\D\ne\G$, choose $\b\in\G\bs\D$. If there
exists $n\in\N$ such that $\a=n\b\in\D$, then we define $\chi(\b)=
\chi(\a)^{1/n}$, otherwise we set $\chi(\b)=1$. Then $\chi$ extends to a
character on the group generated by $\D$ and $\b$, which contradicts the
maximality of $\D$.
\par
First suppose $a_1=a_3=0$. We shall verify that
$$
\psi: x^{\b,\ii}
\mapsto a_4^{-1}\chi(\b)
y^{\b'}(t'_1)^{i_1}(a_2 t'_2)^{i_2}(t'_3)^{i_3}(a_4 t'_4)^{i_4},
\mbox{ where }\b'=\phi(\b)\mbox{ (cf.~(3.7))},
\eqno(3.10)$$
is an isomorphism $(\AA_4,\odot_1)\cong(\AA'_4,\odot_1)$ (and symmetrically,
it is also an isomorphism $(\AA_4,\odot_2)\cong(\AA'_4,\odot_2)$\,).
Suppose $i_1+j_1,i_2+j_2\ge1$ (if otherwise, the verification is easier).
Using (3.9), since $\phi(\a+\b+\si)=\a'+\b'+\si$ by (3.7), we have
$$
\psi(x^{\a,\ii}\odot_1 x^{\b,{\ssc\,}\jj})=
a_4^{-1}\chi(\a+\b+\si)y^{\a'+\b'+\si}u(\a_1t'_1+i_1)(\b_2a_2t'_2+j_2),
\eqno(3.11)$$
where $u=(t'_1)^{i_1+j_1-1}(a_2t'_2)^{i_2+j_2-1}(t'_3)^{i_3+j_3}
(a_4t'_2)^{i_4+j_4}$, and
$$
\psi(x^{\a,\ii})\odot_1\psi(x^{\b,{\ssc\,}\jj})=
a_4^{-2}\chi(\a)\chi(\b)y^{\a'+\b'+\si}u(\a'_1t'_1+i_1)(\b'_2a_2t'_2+a_2j_2)
\eqno(3.12)$$
where the coefficient $a_2$ before $j_2$ is arisen from
$\ptl'_2((a_2t'_2)^{j_2})$. Thus $\psi(x^{\a,\ii})\odot_1\psi(x^{\b,{\ssc\,}\jj})=
\psi(x^{\a,\ii}\odot_1 x^{\b,{\ssc\,}\jj})$ is equivalent to
$a_2(\a_1t'_1+i_1)(\b_2a_2t'_2+j_2)=(\a'_1t'_1+i_1)(\b'_2a_2t'_2+a_2j_2),$
which is obvious since $\a'_1=\a_1,\b'_2=a_2\b_2$. Thus $\psi$ induces an
isomorphism of $\BB\cong\BB'$.
\par
Suppose $a_1\ne0$ or $a_3\ne0$. We define
$$
\psi(x^{\b,\ii})=a_4^{-1}\chi(\b)y^{\b'}(t'_1)^{i_1}
(a_2 t'_2+a_1t'_1)^{i_2}(t'_3)^{i_3}(a_4 t'_4+a_3t'_3)^{i_4}.
\eqno(3.13)$$
Note that if $a_1\ne0$ then $J_2=\{0\}$ or $J_1,J_2\ne\{0\}$ and if
$a_3\ne0$ then $J_4=\{0\}$ or $J_3,J_4\ne\{0\}$, thus the right-hand side
of (3.13) is in $\BB'$. We claim that $\psi$ is an isomorphism
$\BB\cong\BB'$ (but not necessarily an
isomorphism $(\AA,\odot_q)\rar(\AA',\odot_q)$ for $q=1,2$). If we define
two Lie brackets $[\cdot,\cdot]_q$ by $[u,v]_q=u\odot_q v-v\odot_q u$ for
$q=1,2$, then $[u,v]=[u,v]_1+[u,v]_2$. We want to prove
$$
[\psi(x^{\a,\ii}),\psi(x^{\b,{\ssc\,}\jj})]_q=\psi([x^{\a,\ii},x^{\b,{\ssc\,}\jj}]_q),
\eqno(3.14)$$
for $q=1,2$. Again suppose $i_1+j_1,i_2+j_2\ge1$, and suppose $q=1$ (the
proof for $q=2$ is similar). First we calculate
$\psi(x^{\a,\ii}\odot_1 x^{\b,{\ssc\,}\jj})$, which is the term
$$
a_4^{-1} y^{\a'+\b'+\si}
(t'_1)^{i_1+j_1-1}(a_2t'_2+a_1t'_1)^{i_2+j_2-2}
(t'_3)^{i_3+j_3}(a_4t'_4+a_3t'_3)^{i_4+j_4},
\eqno(3.15)$$
where if $i_2+j_2-2<0$ the corresponding factor does not appear,
multiplied by the term
$$
\chi(\a+\b+\si)(\a_1t'_1+i_1)(a_2t'_2+a_1t'_1)(\b_2(a_2t'_2+a_1t'_1)+j_2).
\eqno(3.16)$$
Similarly, $\psi(x^{\a,\ii})\odot_1\psi(x^{\b,{\ssc\,}\jj})$ is (3.15) multiplied by
$$
a_4^{-1}\chi(\a)\chi(\b)
(\b'_2(a_2t'_2+a_1t'_1)+a_2j_2)((\a'_1t'_1+i_1)(a_2t'_2+a_1t'_1)+a_1i_2t'_1),
\eqno(3.17)$$
where the factor $(\b'_2(a_2t'_2+a_1t'_1)+a_2j_2)$ is arisen from
$\ptl'_2(y^{\b'}(a_2t'_2+a_1t'_1)^{j_2})$ and the last factor is arisen
from $\ptl'_1(y^{\a'}(t'_1)^{i_1}(a_2t'_2+a_1t'_1)^{i_2})$. If we denote by
$D(\a,\b)$ the difference between (3.16) and (3.17), then (3.14) is
equivalent to $D(\a,\b)-D(\b,\a)=0$, which is straightforward to verify.
\par
Assume that we have case (ii). We define a new Lie algebra
$\wt{\BB}=\BB(\wt{\G},\wt{J})$ by \vs{-2pt}taking
$$
\matrix{
\wt\G=\{\wt\a=(\a'_3,\a'_4,\a'_1,\a'_2)\,|\,
\a'=(\a'_1,\a'_2,\a'_3,\a'_4)\in\G'\},
\vs{4pt}\hfill\cr
\wt J=\{\wt\ii=(i'_3,i'_4,i'_1,i'_2)\,|\,\ii'=(i'_1,i'_2,i'_3,i'_4)\in J'\}.
\hfill\cr}
\vs{-2pt}\eqno(3.18)$$
It is straightforward to verify that $\wt\psi:\BB'\rar\wt\BB$ defined by
$$
\wt\psi(y^{\a',\ii'})=z^{-\wt\a-\si,\wt\ii},
\eqno(3.19)$$
is a Lie algebra isomorphism, where the symbol $z$ is the symbol $x$ for the
algebra $\wt\BB$. So, it suffices to prove that $\BB\cong\wt\BB$, but using
definition (3.18), case (ii) for the pair of algebras $(\BB,\BB')$ becomes
case (i) for the pair of algebras $(\BB,\wt\BB)$.
\par
``$\Rar$'':
Suppose $\psi:\BB\rar\BB'$ is an isomorphism. Then $\psi$ induces an
isomorphism $\psi:\der\BB\rar\der\BB'$, which maps
$(\der\BB)_{\rm n},(\der\BB)_{\rm s},(\der\BB)_{\rm f}$ respectively to
$(\der\BB')_{\rm n},(\der\BB')_{\rm s}$, $(\der\BB')_{\rm f}$. Denote by
$\BB_{\rm s}$, $\BB_{\rm f}$ the set of {\it ad}-semi-simple,
{\it ad}-locally finite elements of $\BB$ respectively, then by Lemma 3.1,
$\BB_{\rm f}=\F\!+\!\F x^{-\si}$. Thus we have
$$
\pmatrix{
\psi(-x^{-\si})\cr\psi(1)\cr}=G\pmatrix{-y^{-\si}\cr1\cr}
\mbox{ \ for some \ }G=\pmatrix{b_1&b_2\cr b_3&b_4\cr}\in GL_2,
\eqno(3.20)$$
where in general, $GL_n$ is the group of $n\times n$ invertible matrices
over $\F$ and where, we put minus sign before $x^{-\si}$ for later
convenient use. Note that
$$
\{u\in\BB\,|\,[\BB_{\rm f},u]=0\}=
{\rm span}\{x^{\a,\ii}\,|\,\a_2=\a_4=i_2=i_4=0\mbox{ \ or  \ }
\a=\si_1,|\ii|=1\},
\eqno(3.21)$$
where $|\ii|=\sum_{p=1}^4 i_p$ is the {\it level} of $\ii$. Denote the
subalgebra in (3.21) by $\BB^*$. Then we have $\psi(\BB^*)=\BB'^*.$ Observe
that $\BB_{\rm f}=\BB_{\rm s}\Rla(J_2,J_4)=(\{0\},\{0\})$ and that $\BB^*$
is not abelian $\Rla(J_2,J_4)\ne(\{0\},\{0\})$ and $\si\in\G$. Thus we
\vs{-2pt}obtain
$$
\matrix{
(J_2,J_4)=(\{0\},\{0\})\Rla(J'_2,J'_4)=(\{0\},\{0\}),\mbox{ \ and}
\vs{4pt}\hfill\cr
\si\in\G\Rla\si\in\G'\mbox{\ if \ }(J_2,J_4)\ne(\{0\},\{0\}).
\hfill\cr}
\vs{-2pt}\eqno(3.22)$$
If $(J_2,J_4)\ne(\{0\},\{0\})$ and $\si\in\G$, we redenote the {\it derived
subalgebra} $[\BB^*,\BB^*]$ of $\BB^*$ by $\BB^*$. Thus in any case,
$$
\BB^*={\rm span}\{x^{\a,\ii}\,|\,\a_2=\a_4=i_2=i_4=0\},
\mbox{ \ and we have \ }\psi(\BB^*)=\BB'^*.
\eqno(3.23)$$
Note that $J=\{0\}\Rla(\der\BB)_{\rm n}=\{0\}$ and note that
$$
\{u\in\BB\,|\,\ptl(u)=0\mbox{ \ for all \ }\ptl\in(\der\BB)_{\rm n}\}
={\rm span}\{x^{\a,\ii}\,|\,\ii=0\mbox{ \ or \ }\a=\si_1,|\ii|=1\}.
\eqno(3.24)$$
If we denote the derived subalgebra of (3.24) by $\BB^{(0)}$, then
$$
\BB^{(0)}={\rm span}\{x^\a\,|\,\a\in\G\bs\{\si_1,\si_2\}\},
\eqno(3.25)$$
is the Lie algebra $\BB(\G,\{0\})$ (cf.~Theorem 1.1), and we have
$$
J=\{0\}\Rla J'=\{0\},\mbox{ \ and \ }\psi(\BB^{(0)})=\BB'^{(0)}.
\eqno(3.26)$$
We shall study a feature of $\BB^{(0)}$ analogous to that of the Lie
algebras introduced by Zhao [Z]. To do this, take
$\G_{24}=\{(0,\a_2,0,\a_4)\,|\,(\a_1,\a_2,\a_3,\a_4)\in\G\}$ to be the image
of $\G$ under the natural projection $\pi_{24}:\G\rar\F1_{[2]}+\F1_{[4]}$.
Take $\G_{13}=\{\a\in\G\,|\,\pi_{24}(\a)=0\}$, the kernel of $\pi_{24}$.
Then by \vs{-3pt}(2.1),
$$
\BB^{(0)}=\bigoplus_{\l\in\G_{24}}\BB^{(0)}_\l
\mbox{ \ is \,$\G_{24}$-graded with \ }
\BB^{(0)}_\l={\rm span}\{x^\a\,|\,\pi_{24}(\a)=\l\},
\vs{-3pt}\eqno(3.27)$$
and $\BB^{(0)}_0=\BB^{(0)}\cap\BB^*$ and $\psi(\BB^{(0)}_0)=\BB'{}^{(0)}_0$
by (3.23), (3.26), and there is a bijection $\l\mapsto\l'$ from
$\G_{24}\rar\G'_{24}$ such that $\psi(\BB^{(0)}_\l)=\BB'{}^{(0)}_\l$ since
$(\cup_{\l\in\G_{24}}\BB^{(0)}_\l)\bs\{0\}$ is the set of common
eigenvectors of $\ad_{x^{-\si}}$ and $\ad_1$.
\par
{\bf Claim 1}. If $\l\ne0$, then $\BB^{(0)}_\l$ is a cyclic
$\BB^{(0)}_0$-module (i.e., generated by one element), the nonzero scalar
multiples of $x^\a$ for all $\a\in\G\bs\{\si_1,\si_2\}$ with
$\pi_{24}(\a)=\l$ are the only generators.
\par
For any $u=\sum_{i=1}^m c_i x^{\mu_i}$, where $m>1,c_1,...,c_m\in\F\bs
\{0\},\,\pi_{24}(\mu_i)=\l$ and $\mu_1,...,\mu_m$ are distinct, it is
straightforward to verify that
${\rm span}\{\sum_{i=1}^m c_i x^{\a+\mu_i}\,|\,\a\in\G_{13}\}$ is a proper
$\BB^{(0)}_0$-submodule of $\BB^{(0)}_\l$, and that it contains $\la u\ra$,
where in general, we use $\la u\ra$ to denote the cyclic
$\BB^{(0)}_0$-submodule generated by $u$. Thus $u$ is not a generator of
$\BB^{(0)}_\l$ as a $\BB^{(0)}_0$-module.
\par
Let $\a\in\G\bs\{\si_1,\si_2\}$ with $\pi_{24}(\a)=\l$ and $\l\ne0$. Say
$\l_2\ne0$ (the proof for $\l_4\ne0$ is similar). Choose
$\eta=(e,0,-1,0)\in\G,e\ne0,\pm1$ as in the paragraph before (2.47), then
$\eta_k=k\eta+(k-1)\si=(ke+k-1,0,-1,0)\in\G\bs\{\si_1,\si_2\}$ for $k\in\Z$
(except possibly one $k$), and we \vs{-2pt}have
$$
\matrix{
[x^{\eta_k},x^\a]=(ke+k-1)\a_2 x^{\a+k\eta+k\si}\in\la x^\a\ra\mbox{ \ and}
\vs{4pt}\hfill\cr
[x^{\b-k\eta-(k+1)\si},x^{\a+k\eta+k\si}]
=(\b_1-ke-k-1)\a_2 x^{\a+\b}+\b_3\a_4x^{\a+\b-\si}\in\la x^\a\ra,
\hfill\cr}
\vs{-2pt}\eqno(3.28)$$
for all $k\in\Z$ (except possibly two $k$'s) and all $\b\in\G_{13}$. This
proves that
${\rm span}\{x^{\a+\b}\,|\,\b\in\G_{13}\}$ is contained in $\la x^\a\ra$.
Since $\BB_\l^{(0)}={\rm span}\{x^{\a+\b}\,|\,\b\in\G_{13}\}$, we obtain
that $x^\a$ is a generator of $\BB^{(0)}_\l$ as a $\BB^{(0)}_0$-module.
\par
By Claim 1, there exists a bijection $\phi:\a\mapsto\a'$ from
$\G\bs\G_{13}\rar\G'\bs\G'_{13}$ such \vs{-3pt}that
$$
\psi(x^\a)=c_\a y^{\a'}
\mbox{ \ for \ }\a\in\G\bs\G_{13}\mbox{ \ and some \ }c_\a\in\F\bs\{0\}.
\vs{-3pt}\eqno(3.29)$$
Using this and (3.20), we \vs{-3pt}have
$$
c_\a\pmatrix{\a_2\cr a_4\cr}y^{\a'}
=c_\a G\pmatrix{[-y^{-\si},y^{\a'}]\cr[1,y^{\a'}]\cr}
=c_\a G\pmatrix{\a'_2\cr\a'_4\cr}y^{\a'},\mbox{ \ i.e., \ }
\pmatrix{\a_2\cr\a_4\cr}=G\pmatrix{\a'_2\cr a'_4\cr},
\eqno(3.30)$$
for all $\a\in\G\bs\G_{13}$. Let $\tau_2,\tau_4$ be as in (2.2), and write
$\phi(\tau_p)\!=\!\tau'_p\!=\!(\tau'_{p1},\tau'_{p2},\tau'_{p3},\tau'_{p4})$
for $p=2,4$. Then applying $\psi$ to $\a_1ax^{\a+\tau_2+\si}=
[x^\a,x^{\tau_2}],$ we \vs{-2pt}obtain
$$
\matrix{
\a_1ac_{\a+\tau_2+\si}y^{\phi(\a+\tau_2+\si)}
=\!\!\!\!&
c_\a c_{\tau_2}(\a'_1\tau'_{22}-\tau'_{21}\a'_2)y^{\a'+\tau'_2+\si}
\vs{4pt}\hfill\cr&
+c_\a c_{\tau_2}((\a'_3+1)\tau'_{24}-(\tau'_{23}+1)\a'_4)y^{\a'+\tau'_2},
\hfill\cr
}
\vs{-2pt}\eqno(3.31)$$
for $\a,\a+\tau_2\in\G\bs\G_{1,3}.$ The left-hand side has only one term,
thus one of two terms in the right-hand side must be zero for all such $\a$.
This shows that $\tau'_{22}=\tau'_{21}=0$ or $\tau'_{24}=\tau'_{23}+1=0$.
If necessary, by considering the isomorphism
$\wt\phi\cdot\psi:\BB\cong\wt\BB$ instead of $\psi$, where $\wt\BB,\wt\psi$
are defined in (3.18), (3.19), we can suppose
$$
\tau'_{24}=\tau'_{23}+1=0.
\eqno(3.32)$$
By this and (3.30), we have $\pmatrix{a\cr0\cr}=G\pmatrix{\tau'_{22}\cr0\cr}
=\pmatrix{b_1\tau'_{22}\cr b_3\tau'_{22}\cr}$, to give that $b_3=0$.
Similarly, from $(\a_3+1)b x^{\a+\tau_4}=[x^\a,x^{\tau_4}],$ we
\vs{-2pt}obtain
$$
\matrix{
(\a_3+1)bc_{\a+\tau_4}y^{\phi(\a+\tau_4)}
=\!\!\!\!&
c_\a c_{\tau_4}((\a'_1\tau'_{42}-\tau'_{41}\a'_2)y^{\a'+\tau'_4+\si}
\vs{4pt}\hfill\cr&
+((\a'_3+1)\tau'_{44}-(\tau'_{43}+1)\a'_4)y^{\a'+\tau'_4},
\hfill\cr
}
\vs{-2pt}\eqno(3.33)$$
for $\a,\a+\tau_4\in\G\bs\G_{1,3}$, so $\tau'_{42}=\tau'_{41}=0$ or
$\tau'_{44}=\tau'_{43}+1=0$. But from
$\pmatrix{0\cr b\cr}=G\pmatrix{\tau'_{42}\cr\tau'_{44}\cr}$,
we see that $\tau'_{44}\ne0$, \vs{-6pt}thus
$$
\tau'_{42}=\tau'_{41}=0\mbox{ \ and so \ }b_2=0\mbox{ \ and \ }
G=\pmatrix{b_1&0\cr0&b_4\cr}.
\vs{-6pt}\eqno(3.34)$$
Suppose $\a,\b,\a+\b\in\G\bs\G_{13}$, applying $\psi$ to (2.1) with
$\ii=\jj=0$, using (3.30), we \vs{-2pt}have
$$
\matrix{
(\a_1\b_2-\b_1\a_2)c_{\a+\b+\si} y^{\phi(\a+\b+\si)}
+((\a_3+1)\b_4-(\b_3+1)\a_4)c_{\a+\b}y^{\phi(\a+\b)}
\vs{4pt}\hfill\cr
-(\a'_1\b_2-\b'_1\a_2)b_1^{-1}c_\a c_\b y^{\a'+\b'+\si}
-((\a'_3+1)\b_4-(\b'_3+1)\a_4)b_4^{-1}c_\a c_\b y^{\a'+\b'}=0.
\hfill\cr}
\vs{-2pt}\eqno(3.35)$$
This in particular implies that $\phi(\a+\b)=\a'+\b'$ or $\a'+\b'+\si$
for ALL $\a,\b,\a+\b\in\G\bs\G_{13}$ such that
$(\a_3+1)\b_4-(\b_3+1)\a_4\ne0$. From this one can deduce that
$\phi(\a+\b)=\a'+\b'$ for $\a,\b,\a+\b\in\G\bs\G_{13}$ as follows: (3.33)
shows that $\phi(\a+\tau_4)=\a'+\tau'_4$, setting $\b=\a+\tau_4$ in (3.35)
shows that $\phi(2\a)=2\a'$ by noting that the first, third terms are
vanishing; write $\phi(\a+\b)=\a'+\b'+k_{\a,\b}\si$, where $k_{\a,\b}=0,1$;
if $\phi(\a+\b)=\a'+\b'+\si$ for some $\a,\b$, then
$2\a'\!+\!2\b'\!+\!k_{2\a,2\b}\si\!=
\!\phi(2\a)\!+\!\phi(2\b)\!+\!k_{2\a,2\b}\si
\!=\!\phi(2\a\!+\!2\b)\!=\!\phi((\a\!+\!\b)\!+\!(\a\!+\!\b))
\!=\!2(\a'\!+\!\b'\!+\!\si)\!+\!k_{\a+\b,\a+\b}\si$, and
$k_{2\a,2\b}\!=\!2\!+\!k_{\a+\b,\a+\b}\!\ge\!2$, a contradiction. Thus,
$\phi$ can be uniquely extended to a group isomorphism $\phi:\G\rar\G'$.
Using this and applying $\psi$ to $[x^\a,x^{-\a-\si}]=\a_2-\a_4 x^{-\si}$,
by (3.20), (3.34), we obtain
$\a_2b_4\!-\!\a_4b_1y^{-\si}\!=\!c_\a c_{-\a-\si}(\a'_2\!-\!\a'_4y^{-\si})
\!=\!c_\a c_{-\a-\si}(\a_2b_1^{-1}\!-\!\a_4b_4^{-1}y^{-\si})$,
\vs{-5pt}
i.e.,
$$
c_a c_{-\a-\si}=b_1b_4\ \ \for\ \ \a\in\G\bs\G_{13}.
\vs{-5pt}
\eqno(3.36)$$
Comparing the coefficients of $y^{\a'+\b'},y^{\a'+\b'+\si}$ in (3.35), we
\vs{-2pt}obtain
$$
\matrix{
\hfill((\a_3+1)\b_4-(\b_3+1)\a_4)c_{\a+\b}\!\!\!\!&=&\!\!\!\!
((\a'_3+1)\b_4-(\b'_3+1)\a_4)b_4^{-1}c_\a c_\b,
\vs{4pt}\hfill\cr
\hfill(\a_1\b_2-\b_1\a_2)c_{\a+\b+\si}\!\!\!\!&=&\!\!\!\!
(\a'_1\b_2-\b'_1\a_2)b_1^{-1}c_\a c_\b,
\vs{4pt}\hfill\cr
\hfill(\a_3\b_4-\b_3\a_4)c_{-\a-\b-2\si}\!\!\!\!&=&\!\!\!\!
(\a'_3\b_4-\b'_3\a_4)b_4^{-1}c_{-\a-\si}c_{-\b-\si},
\vs{4pt}\hfill\cr
\hfill((\a_1+1)\b_2-(\b_1+1)\a_2)c_{-\a-\b-\si}\!\!\!\!&=&\!\!\!\!
((\a'_1+1)\b_2-(\b'_1+1)\a_2)b_1^{-1}c_{-\a-\si} c_{-\b-\si},
\hfill\cr}
\eqno
\matrix{(3.37)\vs{4pt}\cr(3.38)\vs{4pt}\cr(3.39)\vs{4pt}\cr(3.40)\cr}
\vs{-2pt}$$
for $\a,\b,\a+\b\in\G\bs\G_{13}$, where the last two equations are obtained
from the first two by replacing $\a,\b$ by $-\a-\si,-\b-\si$ respectively.
Setting $\a=\tau_2=(0,a,-1,0)$ in (3.38), (3.39) and multiplying the
obtained two equations, using (3.36), (3.32), we \vs{-2pt}obtain
$$
\matrix{
\hfill\b_1a\b_4\!\!\!&=&\!\!\!-(\tau'_{21}\b_2-\b'_1a)\b_4,
\vs{4pt}\hfill\cr
\hfill-(\b_3+1)b\b_2\!\!\!&=&\!\!\!((\tau'_{43}+1)\b_4-(\b'_3+1)b)\b_2,
\hfill\cr}
\vs{-2pt}\eqno(3.41)$$
for all $\b,\b+\tau_2,\b+\tau_4\in\G\bs\G_{13}$, where the second equation
is obtained analogously from (3.37), (3.40) and by setting
$\a=\tau_4=(0,0,-1,b)$ and using (3.34). Now if we \vs{-5pt}take
$$
a_1=a^{-1}\tau'_{21},\ a_2=b_1^{-1},\
a_3=b^{-1}(\tau'_{43}+1),\
a_4=b_4^{-1},
\vs{-5pt}
\eqno(3.42)$$
(note that $\tau'_{21},\tau'_{43}$ are fixed number since $\tau_2,\tau_4$
are fixed), then by (3.41), (3.30), (3.34), we see that (3.7) holds for
$\b,\b+\tau_2,\b+\tau_4\in\G\bs\G_{13}$ and $\b_2\ne0\ne\b_4$. Since $\phi$
is a group isomorphism, we obtain that (3.7) holds for all $\b\!\in\!\G$.
\vs{-2pt}This proves the theorem if $J\!=\!\{0\}$.
\par
Assume that $J\ne\{0\}$. Since $\psi(x^{-\si})\in\F y^{-\si}$ and
$\ad_{x^{-\si}}$ is semi-simple $\Rla J_2=\{0\}$, thus we obtain that
$J_2=J'_2$. Similarly, $J_4=J'_4$. Assume that $J_1\ne\{0\}=J'_1$. Then we
can find $t_1=x^{0,1_{[1]}}\in\BB^*$ (cf.~(3.23)) with
$[t_1,x^{\tau_2}]=ax^{\tau_2+\si}$, but we can not find $u'\in\BB'^*$
with $[u',y^{\tau'_2}]=y^{\tau'_2+\si}$ because for any
$y^{\a',\ii'}\in\BB'^*$, $[y^{\a',\ii'},y^{\tau'_2}]=
\a'_1\tau'_{22}y^{\tau'_2+\a'+\si,\ii'}$ can not produce a nonzero term
$y^{\tau'_2+\si}$. This is a contradiction. Thus $J_1=J'_1$. Similarly
$J_3=J'_3$. This proves $J=J'$.
\vs{-2pt}\par
It remains to prove that $a_1=0$ if $J_1=\{0\}\ne J_2$ (the proof for
$a_3=0$ if $J_3=\{0\}\ne J_4$ is similar). Consider
$\psi(x^{\tau_2,1_{[2]}})$. Since $\ptl_p(x^{\tau_2,1_{[2]}})
=[u,x^{\tau_2,1_{[2]}}]=0$ for $p=1,4,\,u=1,x^{\tau_2},x^{\tau_4}$ and
$\ptl_q(x^{\tau_2,1_{[2]}})$ $({\rm mod\sc\,}\F x^{\tau_2})$ is equal to a
scalar multiple of $x^{\tau_2,1_{[2]}}$ for $q=2,3$, we obtain that
$\psi(x^{\tau_2,1_{[2]}})\in\BB'^{(0)}+\F y^{\tau'_2,1_{[2]}}$. Write
$\psi(x^{\tau_2,1_{[2]}})=u'+ a'_0 y^{\tau'_2,1_{[2]}}$ for some
$u'\in\BB'^{(0)}$ and $a'_0\in\F\bs\{0\}$, then
$0=\psi([x^{\tau_2},x^{\tau_2,1_{[2]}}])
=c_{\tau_2}([y^{\tau'_2},u']+a'_0\tau'_{21}y^{2\tau'_2+\si})$, but there
does not exists $u'\in\BB'^{(0)}$ such that
$[u',y^{\tau'_2}]=y^{2\tau'_2+\si}$. Thus $\tau'_{21}=0$ and (3.42) gives
that $a_1=0$. This completes the proof
of the theorem.
\vs{-2pt}\qed\par
For $n>0$, denote by $M_{n\times n}$ the algebra of $n\times n$ matrices
with entries in $\F$ and by $GL_n$ the group of invertible $n\times n$
matrices with entries in $\F$. Let
$M=\{\vec m=(m_1,m_2,m_3,m_4)\,|\,m_p=0,1$ for $p=1,2,3,4\}$, corresponding
to 16 possible choices of $J$. For
$\vec m\in M$, let $G_{\vec m}$ be the subgroup of $GL_4$ generated by all
matrices $(^A_C\ ^B_D)$ satisfying
$$
\matrix{
\mbox{either }B=C=0\mbox{ and }A,D\mbox{ have the form }
(^1_{a_1}\ ^0_{a_2}),(^1_{a_3}\ ^0_{a_4}),\mbox{ or else }
\vs{4pt}\hfill\cr
(m_1,m_2)=(m_3,m_4),A=D=0\mbox{ and }B,C\mbox{ have the form }
(^1_{a_1}\ ^0_{a_2}),(^1_{a_3}\ ^0_{a_4}),
\hfill\cr}
\eqno(3.43)$$
such that $a_2,a_4\ne0$ and $a_1=0$ if $m_1=0\ne m_2$ and $a_3=0$ if
$m_3=0\ne m_4.$
\vs{-2pt}\par
Define an action of $G_{\vec m}$ on $\F^4$ by $g(\a)=\a g^{-1}$ for
$\a\in\F^4,g\in G_{\vec m}.$ For any additive subgroup $\G$ of $\F^4$ and
$g\in G_{\vec m}$, the set $g(\G)=\{g(\a)\,|\,\a\in \G\}$ also forms an
additive subgroup of $\F^4$. Denote by $\Omega$ the set of subgroups $\G$ of
$\F^4$ satisfying (1.1). Denote
$\Omega_{\vec m}=\{g(\G)\in\Omega\,|\,g\in G_{\vec m}\}$ and set
$\ol\Omega_{\vec m}=\Omega/\Omega_{\vec m}$, the quotient set of $\Omega$ by
$\Omega_{\vec m}$. Then Theorem 3.2 implies the following theorem.
\vs{-2pt}\par\ni
{\bf Theorem 3.3}. {\it
There exists a 1-1 correspondence between the set of the isomorphism classes
of the simple Lie algebras $\BB(\G,J)$ of Xu type and the following
\vs{-3pt}set
$$
{\cal M}=\{(m_1,m_2,m_3,m_4,\omega)\,|\, \vec
m=(m_1,m_2,m_3,m_4)\in M,\,\omega\in\ol\Omega_{\vec m}\}.
\vs{-3pt}\eqno(3.44)$$ In other word, ${\cal M}$ is the structure
space of the simple Lie algebras $\BB(\G,J)$ of Xu type.}
\qed\par\ \vs{-19pt}\par\ni {\bf
References}\baselineskip=16pt\parskip=5.5pt
\small\par\ni\hi4ex\ha1 [B] Block, R.: On torsion-free abelian
groups and Lie algebras,
 Proc.~Amer.~Math.~Soc. {\bf 9}, 613-620 (1958)
 \par\ni\hi4ex\ha1
[DZ] Dokovic, D. and Zhao, K.: Derivations, isomorphisms and
 second cohomology of generalized Block algebras, Alg.~Colloq. {\bf 3},
 245-272 (1996)
 \par\ni\hi4ex\ha1
[SX] Su, Y. and Xu, X.: Structure of divergence-free Lie algebras,
J.Alg. {\bf243} (2001), 557-595.
 \par\ni\hi4ex\ha1
[SZ] Su, Y. and Zhou, J.: Structure of the Lie algebras with a feature of
 Block algebras, Comm.~Alg. {\bf30} (2002), 3205-3226.
 \par\ni\hi4ex\ha1
[X] Xu, X.: Generalizations of Block algebras, Manuscripta Math. {\bf 100},
 489-518 (1999)
 \par\ni\hi4ex\ha1
[Z] Zhao, K.: A Class of infinite dimensional simple Lie algebras,
 J.~London Math.~Soc. (2) {\bf62}, 71-84 (2000)
\end{document}